\documentclass[11pt,a4paper]{amsart}
\usepackage{amsmath,amscd,amsfonts,amssymb,mathrsfs,amsthm,bbm,enumerate,psfig}
\usepackage{hyperref}
\usepackage[latin1]{inputenc}
\newtheorem{prop}{Proposition}[section]
\newtheorem{thm}[prop]{Theorem}
\newtheorem{lem}[prop]{Lemma}
\newtheorem{cor}[prop]{Corollary}
\newtheorem{defi}[prop]{Definition}
\theoremstyle{definition}
\newtheorem{rem}[prop]{Remark}

\newcommand{\D}{\mathbb D}
\newcommand{\R}{\mathbb R}
\newcommand{\C}{\mathbb C}
\newcommand{\N}{\mathbb N}
\newcommand{\Z}{\mathbb Z}
\renewcommand{\S}{\mathbb S}

\newcommand{\bp}{b^+_2}
\newcommand{\set}[1]{\{#1\}}
\newcommand{\suchthat}{\;|\;}
\newcommand{\restrict}{\,|} 
\newcommand{\Poincare}{Poincar\'{e}\;}
\newcommand{\define}{\;{\rm :=}\;}

\newcommand{\mt}[1]{\text{\rm #1}}
\newcommand{\mtq}[1]{\qquad\text{\rm #1}\;\;}
\newcommand{\mtqq}[1]{\qquad\text{\rm #1}\qquad}

\newcommand{\pd}{\mt{PD}}
\newcommand{\ie}{i.e.,\ }
\newcommand{\cf}{cf.\ }
\newcommand{\lra}{\longrightarrow}

\newcommand{\Niederkruger}{Niederkr\"uger\;}
\newcommand{\Ozbagci}{\"Ozba\u{g}c\i\;}
\newcommand{\Ozturk}{\"Ozt\"urk\;}
\newcommand{\Schonenberger}{Sch\"onenberger\;}
\DeclareMathOperator{\Sl}{\mathrm{SL}}
\DeclareMathOperator{\uni}{\mathrm{U}}
\DeclareMathOperator{\psl}{\mathrm{PSL}}
\DeclareMathOperator{\inte}{\mathrm{Int}}
\DeclareMathOperator{\rot}{\mathtt{rot}}
\DeclareMathOperator{\tb}{\mathtt{tb}}
\DeclareMathOperator{\lk}{\mathtt{lk}}
\DeclareMathOperator{\diag}{\mathrm{diag}}

\DeclareMathOperator{\pc}{\mathrm{PC}}
\DeclareMathOperator{\cav}{\mathrm{cav}}
\DeclareMathOperator{\odd}{\mathrm{odd}}

\begin{document} 

\author{Kai Zehmisch}

\address{Mathematisches Institut, Universit\"at Leipzig, Augustusplatz 10-11, 
D-04109, Leipzig, Germany}

\email{Kai.Zehmisch@math.uni-leipzig.de}

\title{Strong fillability and the Weinstein conjecture}

\date{April 8, 2005}

\begin{abstract}
  Extending work of Chen, we prove the Weinstein conjecture in dimension three for strongly fillable
  contact structures with either non-vanishing first Chern class or with strong and exact filling
  having non-trivial canonical bundle.
  This implies the Weinstein conjecture for certain Stein fillable contact structures obtained by the 
  Eliashberg-Gompf construction.
  For example we prove the Weinstein conjecture for the Brieskorn homology spheres $\Sigma(2,3,6n-1)$,
  $n\geq2$, oriented as the boundary of the corresponding Milnor fibre.
  Furthermore, for tight contact structures on odd lens spaces, non-contractible closed Reeb orbits are found.
\end{abstract}



\maketitle

\tableofcontents


\section{Basic definitions and main results}

In order to formulate our results, we need to make precise certain well-known notions.
\smallskip\\
{\bf Co-orientable contact $3$-manifolds.}
Let $(M,\xi)$ be a co-orientable contact $3$-manifold, this is a $3$-manifold $M$ equipped
with a co-orientable contact structure $\xi$; the latter is by definition the kernel of
some $1$-form $\lambda$ on $M$ such that the $3$-form $\lambda\wedge d\lambda$ is a
volume form.
We call each such $\lambda$ a \emph{$\xi$-defining contact form}\/.
We equip $M$ with the orientation induced by $\xi$, \ie with the orientation defined by the
non-vanishing top-rank form $\lambda\wedge d\lambda$, where $\lambda$ is any
$\xi$-defining contact form. Moreover, we equip $\xi$ with a (bundle) orientation and
denote the resulting oriented plane field by $\xi_+$. (Given the orientation of $M$,
defining an orientation of $\xi$ is equivalent to defining a co-orientation of $\xi$,
\ie an orientation of the line bundle $TM/\xi$.) There exists a $\xi$-defining contact
form $\lambda$ such that the top-rank form $d\lambda|_{\xi}$ induces the given
orientation on $\xi$; we call each such $\lambda$ a \emph{$\xi_+$-defining contact form}\/.
We can identify $\xi_+$ with the {\it positive conformal class}\/ $\pc(\xi_+)$ consisting
of all $\xi_+$-defining contact forms. The group $C^\infty(M,\R_{>0})$ of all
positive-valued functions on $M$ acts freely and transitively on $\pc(\xi_+)$ by
multiplication.
\smallskip\\
{\bf $\lambda$-Reeb orbits and links.}
For each element $\lambda$ of $\pc(\xi_+)$, the $\lambda$-{\it Reeb vector field}\/
$X_{\lambda}$ is defined to be the unique vector field $X$ on $M$ with $\iota_{X}d\lambda=0$
and $\lambda(X)=1$. 
The integral curves of $X_{\lambda}$ are canonically oriented by $X_{\lambda}$, and are called
$\lambda$-{\it Reeb orbits}\/.
A finite, non-empty union of (necessarily embedded and disjoint) closed $\lambda$-Reeb orbits
(\ie $\lambda$-Reeb orbits diffeomorphic to $\S^1$) is called a $\lambda$-{\it Reeb link}\/.
\smallskip\\
{\bf Hypersurfaces of contact type.}
Consider a compact symplectic $4$-manifold $(W,\omega)$ (possibly with boundary).
We orient it by the volume form $\omega\wedge\omega$.
A closed hypersurface $M$ in $W$ is said to be of {\it contact type}\/ if there exists a (necessarily
(co-)oriented) contact structure $\xi_+$ on $M$ and a contact form $\lambda\in\pc(\xi_+)$
satisfying $d\lambda=i^*\omega$, where $i:M\lra W$ denotes the inclusion map.
We write $(M,\lambda)$ for the hypersurface of contact type to indicate a particular choice
of a contact form $\lambda$ and we write $(M,\xi_+)$ if only the co-oriented 
contact structure is used.
\smallskip\\
{\bf Strong fillings.}
If the boundary $\partial W$ of $W$ is non-empty, we orient $\partial W$ by
$\iota_{\nu}(\omega\wedge\omega)$, where $\nu$ is an outward pointing vector field.
If $\partial W=M$ as oriented manifolds we call $M$ {\it strongly convexly fillable}\/ and
$(W,\omega)$ a {\it strong convex filling}\/ of $(M,\lambda)$.
If $\partial W=\overline{M}$, where $\overline{M}$ is the manifold $M$ with reverse orientation,
we call $M$ {\it strongly concavely fillable}\/ and $(W,\omega)$ a {\it strong concave filling}\/ of
$(M,\lambda)$.
We wish to point out that both concepts, being of `contact type' and being a `strong filling', correspond
to the same convexity condition, the difference between the two concepts is only topological
(hypersurface vs. oriented boundary).
\smallskip\\
Denoting by $c_1(\xi_+)\in H^2(M,\Z)$ the first Chern class of the oriented $2$-plane
bundle $\xi_+$ on $M$ (\ie the first Chern class with respect to any complex structure on $\xi_+$
which is compatible with the orientation), we state our first result as:

\begin{thm} 
  \label{Chen1}
  Let $(M,\xi_+)$ be a closed co-oriented strongly convexly fillable contact manifold and $\lambda$
  a $\xi_+$-defining contact form in $\pc(\xi_+)$.
  If $c_1(\xi_+)\neq 0$ then there is a $\lambda$-Reeb link whose homology class is
  \Poincare dual to $-c_1(\xi_+)$.
\end{thm}

Before we come to the second result we shall make the following remark:\;
the first Chern class $c_1(W,\omega)$ of the symplectic manifold $(W,\omega)$ is by definition the first
Chern class of the tangent bundle of $W$ equipped with any $\omega$-compatible almost complex structure.
Fixing an $\omega$-compatible almost complex structure $J$, the {\it canonical bundle}\/ $K$ on $W$
is defined to be $\Lambda^2T^{1,0}$, where $T^*W\otimes\C=T^{1,0}\oplus T^{0,1}$ is the eigenspace
decomposition with respect to to eigenvalues $i$ and $-i$ of the induced action of $J$.
Notice that $c_1(K)=-c_1(W,\omega)$ (see \cite[Example 21.7]{botttu}).

\begin{thm} 
  \label{Chen2}
  Let $(M,\xi_+)$ be a closed co-oriented contact manifold and $(W,\omega)$ a strong convex filling.
  Suppose that
  $c_1(W,\omega)\neq 0$ and $\omega$ is exact.
  Then for any $\xi_+$-defining contact form $\lambda\in\pc(\xi_+)$ there exists a closed
  $\lambda$-Reeb orbit.
\end{thm}

\begin{rem}
  \label{comment}
  In fact, in the situation of Theorem \ref{Chen2}, there is a $\lambda$-Reeb link whose homology
  class is \Poincare dual to $-c_1(\xi_+)$ regardless whether the first Chern class of the contact
  structure vanishes or not.
  If $(W,\omega)$ is a strong convex filling of $(M,\xi_+)$, as in the situation of Theorem \ref{Chen1},
  then the identity $i^*c_1(W,\omega)=c_1(\xi_+)$ holds, where $i:M\lra W$ denotes the inclusion map, because
  of the isomorphism $TW\restrict_M\cong\xi\oplus\underline{\C}$ as complex vector bundles over $M$.
  In particular, if $c_1(\xi_+)\neq0$ then also $c_1(W,\omega)\neq0$, meaning that if additionally
  $\omega$ is exact then Theorem \ref{Chen2} follows from Theorem \ref{Chen1}.
  But it maid be possible that $c_1(\xi_+)=0$ without $c_1(W,\omega)$ being zero.
\end{rem}

Theorems \ref{Chen1} and \ref{Chen2} confirm certain cases of the Weinstein conjecture in the
$3$-dimensional case (\cf \cite{wstrab}); it asks for a closed $\lambda$-Reeb orbit for {\it every}\/ contact form
$\lambda$ on a closed $3$-manifold.
We will say that the {\it Weinstein conjecture holds true}\/ for $M$ (for the oriented manifold $M$, for
{\it the co-orientable contact structure}\/ $\xi$) if there is a closed $\lambda$-Reeb orbit {\it for all}\/ contact
forms $\lambda$ on $M$ ({\it for all}\/ contact forms $\lambda$ on $M$ such that $\lambda\wedge d\lambda$
induces the given orientation on $M$, {\it for all}\/ $\xi_+$-defining contact forms $\lambda\in\pc(\xi_+)$)
respectively.
In this language Theorem \ref{Chen1} and Theorem \ref{Chen2} show that under the posed conditions the
Weinstein conjecture holds true for the contact structure $\xi$.
\smallskip\\
If $(M,\lambda)$ is a hypersurface of contact type in $(W,\omega)$ then one can consider the
characteristic line bundle $\mathcal{L}_M=\ker(i^*\omega)$, where $i:M\lra W$ denotes
the inclusion map.
In this case the Weinstein conjecture asks for a closed characteristic of $\mathcal{L}_M$.
For more about the extrinsic form of the conjecture and the state of the art of this problem we refer
to \cite{hz,schlfrau,ginz}.
\smallskip\\
In \cite{h2} Hofer proved the Weinstein conjecture for $\S^3$ (using Rabinowitz' periodic orbit theorem,
see \cite{rab,wstrab}, and his results for overtwisted contact structures in \cite{h2}) and for all closed
$3$-manifolds $M$ with non-trivial second homotopy group.
Notice that any $3$-manifold $B$ which is covered by such an $M$ satisfies the Weinstein conjecture.
For example the Weinstein conjecture holds true for all lens spaces $L_{p,q}$ for $p>q\geq1$ coprime
(notice that $L_{0,1}=\S^1\times\S^2$ and $L_{1,1}=\S^3$) and the \Poincare homology sphere
$\Sigma(2,3,5)$, which are universally covered by $\S^3$.
Further, Hofer \cite{h2} verified the Weinstein conjecture for all overtwisted contact structures $\xi$
on closed $3$-manifolds.
In particular the Weinstein conjecture holds true for all virtually overtwisted contact structures, \ie for all
tight contact structures for which the lift to a finite cover becomes overtwisted.
\smallskip\\
A co-orientable contact structure is called {\it planar} if there exists a supporting open book decomposition of
the underlying closed $3$-manifold (in the sense of Giroux, \cf \cite{planarlec}) which has genus zero pages
(\cf \cite{ach,planar}).
In \cite{ach} Abbas, Cieliebak and Hofer verified the Weinstein conjecture for planar contact structures.
In Section \ref{brie} we will prove:

\begin{thm}
  The Weinstein conjecture holds true for the positively oriented Brieskorn homology spheres
  $+\Sigma(2,3,6n-1)$, $n\geq2$.
\end{thm}

All closed $\lambda$-Reeb orbits found in \cite{h2} are contractible (a covering induces an
injection on $\pi_1$) and the $\lambda$-Reeb links obtained in \cite{ach} are zero in integral
homology of the underlying closed $3$-manifold.
In Section \ref{lens} below we will establish non-contractible (and in fact not null-homologous)
$\lambda$-Reeb orbits for so called odd lens spaces (proving Theorem \ref{nclroools} below).
We call the lens space $L_{p,q}$ {\it odd}\/, and will write $L_{p,q}^{\odd}$, if there is at least one odd
integer $n_i$ in the associated continued fraction expansion $[n_1,\ldots,n_k]$ of the fraction
$-\frac{p}{q}$.

\begin{thm}
  \label{nclroools}
  Let $p>q\geq1$ be coprime integers.
  For all odd lens spaces $L_{p,q}^{\odd}$ and all tight contact structures $\xi$ there exists a
  non-contractible $\lambda$-Reeb orbit for all $\lambda\in\pc(\xi_+)$.
\end{thm}

The proof of Theorem \ref{Chen1} and Theorem \ref{Chen2} has two main ingredients, which we
shall state next.
The first ingredient is due to Chen \cite{weichen}.
He proved the Weinstein conjecture for particular classes of contact type hypersurfaces in
$4$-manifolds using work of Taubes \cite{existence} on Seiberg-Witten equations and
pseudo-holomorphic curves as well as stretching the neck, which is due to Hofer, Wysocki and
Zehnder \cite{s3tight}.
Some ideas in \cite{weichen} are borrowed from \cite{h2,egh}.
\smallskip\\
For a compact $4$-dimensional manifold $W$ (possibly with boundary) denote by $\bp(W)$ the
number of positive eigenvalues of the intersection form $Q_W$ of $W$ (see 
\cite[Definition 1.2.1]{kirby}).

\begin{thm}[Chen]
  \label{base}
  Let $(W,\omega)$ be a closed connected symplectic $4$-manifold with $\bp(W)>1$ and let
  $(M,\lambda)$ be a hypersurface of contact type in $(W,\omega)$.
  Set $\xi\define\ker \lambda$.
            \begin{enumerate}
            \item
              \label{base1}
              If $c_1(\xi_+)\neq 0$ then there exists a $\lambda$-Reeb link whose homology class
              is \Poincare dual to $-c_1(\xi_+)$.
            \item
              \label{base2}
              If $M$ bounds a submanifold $\widehat{W}$ of $W$
              with $c_1(\widehat{W})\neq 0$ and
              $\omega\restrict_{\widehat{W}}$ is exact, then there exists a closed $\lambda$-Reeb orbit.
              In fact there is a $\lambda$-Reeb link whose homology class is \Poincare dual to
              $-c_1(\xi_+)$.
            \end{enumerate}
\end{thm}

Notice that Theorems \ref{Chen1} and \ref{Chen2} remove the assumption $\bp(W)>1$ in this
result and establish its conclusion for all $\xi_+$-defining contact forms $\lambda$; not only for
those with $d\lambda=i^*\omega$, where $i:M\lra W$ denotes the inclusion map
(in the sequel we will also write $\omega\restrict_M$ for $i^*\omega$).
\smallskip\\
The second ingredient is that under certain circumstances a closed co-orientable
contact $3$-manifold $(M,\xi)$ can be realized as a hypersurface of contact type
in a closed symplectic $4$-manifold or in a compact symplectic $4$-manifold with
boundary equal to $(M,\xi_+)$.
The following result is due to Etnyre and Honda \cite[Theorem 1.3]{ehdscob}, see also
\cite{tightsw,explicit,aostein,few,etnfll} for previous work in this direction.

\begin{thm}[Etnyre, Honda]
  \label{concave}
  Any closed connected contact $3$-manifold admits a connected strong concave filling.
  In fact there are infinitely many strong concave fillings which are mutually not homotopy
  equivalent and not related by a sequence of blow-downs and blow-ups.
\end{thm}

This article is organised as follows:\;
in Section \ref{appex} we translate our criterion to decide whether the Weinstein conjecture holds true
for a Stein fillable contact structure given by Theorem \ref{Chen1} or Theorem \ref{Chen2} into the language of
Legendrian $(-1)$-surgery (see Corollary \ref{wcegconstr}).
This leads to a proof of the Weinstein conjecture for the positively oriented Brieskorn homology spheres
$+\Sigma(2,3,6n-1)$, $n\geq2$, (see Section \ref{brie}).
In Section \ref{lens} we treat the odd lens spaces and prove Theorem \ref{nclroools}.
The main Theorems are established in two steps.
The first one is made in Section \ref{dircob} and the second in Section \ref{rpmt}.
\smallskip\\
\textit{Acknowledgements.}
The research presented in this article was carried out while I was supported by the DFG
through the Graduiertenkolleg - Analysis, Geometrie und ihre Verbindung zu den
Naturwissenschaften at the Universit\"at Leipzig.
I am grateful to my supervisor Prof.~M.~Schwarz for his steady help and encouragement.
I would like to thank Paolo Ghiggini and Stephan \Schonenberger for the stimulating e-mail conversation
as well as John Etnyre for explaining me the proof of Theorem \ref{concave}.
I am indebted to Casim Abbas, Peter Albers, Shahram Biglari, Dragomir Dragnev, Frank Klinker, Otto van~Koert,
Matthias Kurzke, Rainer Munck, Marc Nardmann, Klaus \Niederkruger\!\!, Burak \Ozbagci\!\!, Ferit \Ozturk and
Felix Schlenk for many helpful suggestions and comments.


\section{Application to Stein fillable contact structures \label{appex}}

In this section we describe the Eliashberg-Gompf construction of Stein fillable contact structures.
It provides us with examples having computable topological invariants useful to decide whether
the constructed contact structure satisfies the Weinstein conjecture.
\smallskip\\
A {\it Stein $4$-manifold with boundary}\/ (or simply a {\it Stein surface with boundary}\/) is a triple
$(W,J,\varphi)$ consisting of a smooth $4$-manifold $W$ with non-empty boundary,
a complex structure $J$ on $W$ such that there exists $N\in\N$ with the property that
$(\inte(W),J)$ is biholomorphically equivalent to a complex submanifold of $\C^N$, and a
Morse function $\varphi:W\longrightarrow \R$ such that $\varphi\restrict_{\partial W}$
is constant and the $2$-form $\omega_{\varphi}=-dJ^*d\varphi$ defines a symplectic structure
on $W$, where $J^*\alpha =\alpha\circ J$ for all $1$-forms $\alpha$ on $W$.
For any non-empty regular level set $\varphi^{-1}(c)$ the $1$-form
$\lambda =-J^*d\varphi\restrict_{\varphi^{-1}(c)}$ defines a co-oriented contact structure
$\xi_{J}$ on $\varphi^{-1}(c)$.
If the strong convex filling $(W,\omega)$ of $(M,\lambda)$ {\it carries a Stein structure}\/, that is,
there exists a  Stein $4$-manifold with boundary $(W,J,\varphi)$ such that
$\omega=\omega_{\varphi}$ and $\lambda=-J^*d\varphi\restrict_{\partial W}$,
we call $(W,\omega)$ a {\it Stein filling}\/ of $(M,\lambda)$.
It is unknown whether strong convex fillability implies Stein fillability of connected contact 
$3$-manifolds.
But there are examples of disconnected strongly convexly fillable contact manifolds,
such that their corresponding fillings cannot carry any Stein structure (see \cite{mcduconty,exageig}).
\smallskip\\
Suppose that the Stein surface with boundary $(W,J,\varphi)$ is a handlebody with only one
$0$-handle and $m$ $1$-handles where $m\geq 0$.
The induced contact structure $(\partial W, \xi_{J})$ is contactomorphic to
$(\#m\S^1\times\S^2,\xi_0)$ where $\xi_0$ denotes the standard contact structure on
$\#m\S^1\times\S^2$, and  $(W,J,\varphi)$ is the unique Stein filling of $(\#m\S^1\times\S^2,\xi_0)$
(\cite{elfll,el1}).
\smallskip\\
In the remainder of this section we assume that the reader is familiar with \cite{gompf}.
We call a link $L=(K_1,\ldots,K_n)$, $n\in\N$, in a contact 3-manifold $(M,\xi)$ {\it Legendrian}\/  if 
the knots $K_i$ are tangent to $\xi$.
Any Legendrian link in $(\#m\S^1\times\S^2,\xi_0)$ is contact isotopic to a Legendrian link in
so called {\it standard form}\/ (see \cite[Definition 2.1 and Theorem 2.2]{gompf}).
If a Legendrian knot $K$ is in standard form one can define its  {\it Thurston-Bennequin invariant}\/
$\tb(K)$.
In the special case that $K$ is null-homologous, $\tb(K)$ is the linking-number of the knot $K$
and the parallel push-off knot determined by the {\it canonical framing}\/, this is any vector
field along $K$ transverse to $\xi_0$ respecting the co-orientation.
In exactly the same way any normal vector field to $K$ defines a {\it framing}\/.
\smallskip\\
With a Legendrian link $L$ in standard form one can associate a second invariant -- the
{\it rotation number}\/ $\rot(L)$  -- defined in \cite[Section 2]{gompf} or \cite[Formula 1.2]{gompf}.
In the case of a Legendrian link in $(\S^3,\xi_0)$ the rotation number $\rot(L)$ equals the relative
Chern number $\langle c_1(\xi_0,\tau),[F]\rangle$ of $\xi_0$ relative to a tangent vector field $\tau$
along $L$, evaluated on a Seifert surface $F$ (\ie $\partial F =L$ and
$c_1(\xi_0,\tau)\in H^2(\S^3,L;\Z)$).
In other words, $\rot(L)$ is the degree of $\tau$ with respect to any trivialisation of $\xi_0\restrict_F$.
(This is independent of the choice of a Seifert surface $F$ because $c_1(\xi_0)$ vanishes.)
The following Theorem is due to Eliashberg \cite{stein,legtranstifgtel} and Gompf \cite{gompf}.

\begin{thm}[Eliashberg, Gompf]
  \label{egconstr}
  Let $W$ be an oriented compact connected $4$-manifold with non-empty boundary.
  Then $W$ admits the structure of a Stein surface with boundary if and only if it carries a
  handlebody decomposition with the following properties:
                    \begin{enumerate}
                    \item 
                      There are no $3$- and $4$-handles.
                    \item 
                      $W$ is built from the unique Stein filling of $(\#m\S^1\times\S^2,\xi_0)$
                      by attaching $2$-handles $h_i$,
                      $i=1,\ldots,n$, to $K_i$ with framing $\tb(K_i)-1$, where $L=(K_1,\ldots,K_n)$, $n\in\N$,
                      is a Legendrian link in $(\#m\S^1\times\S^2,\xi_0)$ in standard form.
                    \end{enumerate}
  The handle decomposition of the Stein structure $(W,J,\varphi)$ is induced by $\varphi$ and
  the first Chern class $c_1(W,\omega_{\varphi})\in H^2(W;\Z)$ is represented by a cocycle whose value
  on $[D_i]$ (the class of the core of $h_i $ in $H_2(W;\Z)$ oriented as at the end of \cite[Section 1]{gompf})
  is equal to $\rot(K_i)$.
\end{thm}

In the situation of Theorem \ref{egconstr} we will say that $(W,J,\varphi)$ is obtained from
(the unique Stein filling of)
$(\#m\S^1\times\S^2,\xi_0)$ by {\it Legendrian $(-1)$-surgery}\/ along $L$.
The surgered contact manifold $(\partial W,\xi_J)$ is Stein cobordant to $(\#m\S^1\times\S^2,\xi_0)$
(see Section \ref{dircob}) and Stein filled by $(W,\omega_{\varphi})$.
Notice that \cite[Theorem 1.2]{tightsw} implies that if $J_1$ and $J_2$ are two Stein structures on
$W$ with $c_1(W,J_1)\neq c_1(W,J_2)$ then the induced contact structures $\xi_{J_1}$ and $\xi_{J_2}$
on $\partial W$ are not isotopic.

\begin{rem}
  \label{compute}
  Using \cite[p.~ 49]{dgs} and \cite[p.~ 658]{gompf} we have in more explicit terms that
                                                            \begin{gather}
                                                              \label{orient1}
                                                              \pd \big(c_1(W,\omega_{\varphi})\big)=
                                                              \sum_{i=1}^n\rot(K_i)[C_i]
                                                              \;,
                                                            \end{gather}
  where $[C_i]$ is the class of the cocore of the $2$-handle $h_i $ in $H_2(W,\partial W;\Z)$
  provided with the orientation mentioned in Theorem \ref{egconstr}.
  By Remark \ref{comment} we also get
                                                            \begin{gather}
                                                              \label{orient2}
                                                              \pd \big(c_1(\xi_J)\big)=
                                                              \sum_{i=1}^n\rot(K_i)[\partial C_i]
                                                              \;,
                                                            \end{gather}
  which is an element of $H_1(\partial W;\Z)$.
  If the surgery is performed on $\S^3$ then the $[C_i]$ freely generate $H_2(W,\partial W;\Z)$ and
  the $[\partial C_i]$ generate $H_1(\partial W;\Z)$ with relations
                                                            \begin{gather}
                                                              \label{lkrel}
                                                              \sum_{j=1}^n\lk(K_i,K_j)[\partial C_j]=
                                                              0
                                                              \;,
                                                            \end{gather}
  where $\big(\lk(K_i,K_j)\big)$ is the linking matrix of the link $L$.
  In this special case (surgery along $\S^3$), the orientation mentioned in Theorem \ref{egconstr} can
  be described as follows:\;
  orienting the knot $K_j$ (for example as the boundary of an oriented Seifert surface $F_j$) together with
  the canonical orientation of the ambient manifold $\S^3$ gives an orientation
  of a small normal disc $N_j$ to $K_j$.
  We choose the orientation of the cocore $C_j$ such that the normal disc $N_j$ will represent
  the class $[C_j]$ with the same sign.
  Analogously, the boundary orientation of $\partial N_j$ determines the sign of $[\partial C_j]$.
  It turns out to be useful to remark that the change of the orientation of $K_j$ induces the
  multiplication of $[C_j]$, $[\partial C_j]$ and $\rot(K_j)$, respectively, by the factor $-1$.
  Therefore, the signs of the summands in both equations \eqref{orient1} and \eqref{orient2} remain
  unchanged.
  But the sign of $\tb(K_j)$ and the diagonal elements in $\big(\lk(K_i,K_j)\big)$ are not effected
  by orientation-reversing $K_j$.
  This means that the signs of the summands in equation \eqref{lkrel} do not change either
  (of course $\lk(K_i,K_j)$ changes sign for $i\neq j$).
\end{rem}

Putting Theorems \ref{Chen2} and \ref{egconstr} together we obtain

\begin{cor}
  \label{wcegconstr}
  Let $(W,J,\varphi)$ be a Stein surface with boundary obtained from the unique Stein filling of
  $(\#m\S^1\times\S^2,\xi_0)$
  via Legendrian $(-1)$-surgery along $L=(K_1,\ldots,K_n)$, $n\in\N$.
  Suppose that the rotation number $\rot(K_i)$ does not vanish for at least one $i\in\{1,\ldots,n\}$.
  Then the Weinstein conjecture holds true for the induced contact structure $\xi_J$ on $\partial W$.
  In fact there exists a $\lambda$-Reeb link whose homology class equals
                                                            \begin{gather*}
                                                              -\sum_{i=1}^n\rot(K_i)[\partial C_i] 
                                                            \end{gather*}
  for all $\lambda\in\pc(\xi_J)$.
\end{cor}


\section{The Weinstein conjecture for the Brieskorn spheres\\
  $\;\;+\Sigma(2,3,6n-1)$, $n\geq2$\label{brie}}

Let $n$ be a natural number and let $\varepsilon\in\set{0<|z|<1}\subset\C$ be a fixed parameter.
Following \cite[p.~74]{kirby} we define the {\it Brieskorn manifold}\/ $\Sigma(2,3,6n-1)$ as the oriented
boundary of the compactified {\it Milnor fibre}
                                                            \begin{gather*}
                                                              M_c(2,3,6n-1)=
                                                              \set{(x,y,z)\in\C^3\suchthat
                                                                x^2+y^3+z^{6n-1}=
                                                                \varepsilon}
                                                              \cap\D^6
                                                              \;.
                                                            \end{gather*}
Alternatively, $-\Sigma(2,3,6n-1)=\overline{\Sigma(2,3,6n-1)}$ is the oriented boundary of the
nucleus $N(n)$ (\cf \cite{nuclei}).
The intersection form $Q_{N(n)}$ is unimodular (see \cite[Figure 8.14]{kirby}).
Therefore, by \cite[Corollary 5.3.12 and Remark 1.2.11]{kirby} the Brieskorn manifold
$\pm\Sigma(2,3,6n-1)$ is in fact a homology sphere  (see Remark \ref{bpeq}) and is therefore
called a {\it Brieskorn homology sphere}\/.
\smallskip\\
If $n\geq2$ (the case $n=1$ corresponds to the \Poincare homology sphere) then
$\pi_1\big(\Sigma(2,3,6n-1)\big)$ is infinite (see \cite[Section 1.1.3]{iohspheres}) and not Abelian.
Because $-\Sigma(2,3,6n-1)$ admits a description as {\it Seifert fibred homology sphere}\/
$M(-\frac12,\frac13,\frac{n}{6n-1})$ (see \cite[Section 1.1.4]{iohspheres} and \cite{3m}) it follows that
$\Sigma(2,3,6n-1)$ is irreducible (see \cite[Proposition 1.12]{3m}) and $\pi_i\big(\Sigma(2,3,6n-1)\big)=0$
for all $i>1$ (see \cite[Corollary 3.9]{3m}).
In particular, $\Sigma(2,3,6n-1)$ cannot be covered by a homotopy sphere.
In fact, it follows from \cite[Section 1]{latt2} that the universal cover of $\Sigma(2,3,6n-1)$ is the
universal cover of $\psl(2,\R)$ and hence equal to $\R^3$ (use that $\Sl(2,\R)$ is diffeomorphic to
$\S^1\times\R^2$).
\smallskip\\
Let $n=2$ and consider $\Sigma(2,3,11)$ (which is equal to $M(\frac12,-\frac13,-\frac{2}{11})$ in
the notation above).
By \cite[Theorem 4.4]{2311} there exist, up to isotopy, exactly two tight contact structures
$\xi_{\pm}$ on $\Sigma(2,3,11)$, both of which are Stein fillable.
The tight contact structures $\xi_{\pm}$ are obtained by Legendrian $(-1)$-surgery along a
Legendrian link $L$ in $\S^3$ with Legendrian knots having $\rot=0$, except for exactly one Legendrian
knot which has $\rot=\pm1$ (see \cite[Section 4.1.4]{2311}).
As Stephan \Schonenberger pointed out to author, for $n\geq3$ a similar statement is true,
as can be seen by using the methods developed in \cite{2311}.
On $\Sigma(2,3,6n-1)$ there exist, up to isotopy, exactly two tight contact structures, both
constructed by Legendrian $(-1)$-surgery along a Legendrian link $L$ in $\S^3$ having $\rot=\pm1$.
Combining Corollary \ref{wcegconstr} with the verified Weinstein conjecture for overtwisted
contact structures confirmed in \cite{h2} and the above mentioned classification theorem (notice
that if the Weinstein conjecture holds true for $\xi$ then the Weinstein conjecture follows trivially
for all contact structures contactomorphic to $\xi$) we have

\begin{cor}
  \label{briecor}
  The Weinstein conjecture holds true for $+\Sigma(2,3,6n-1)$, $n\geq2$.
\end{cor}

By the above discussion this result is not covered by \cite{h2}.
Indeed, as Paolo Ghiggini pointed out to author, the tight contact structures on $+\Sigma(2,3,6n-1)$,
$n\geq2$, are universally tight:\;
by \cite[Theorem 1.3(a) and Corollary 2.2]{transfol} there exists a universally tight contact structure
on $+\Sigma(2,3,6n-1)$, which must be isotopic to either $\xi_-$ or $\xi_+$.
The contact structure $\bar{\xi}_{\mp}$, which is the contact structure $\xi_{\mp}$ with
orientation reversed, is isotopic to $\xi_{\pm}$.
Because this operation, called {\it conjugation}, preserves (universal) tightness, the claim follows.
(In fact, the non-isotopic tight contact structures $\xi_-$ and $\xi_+$ on $+\Sigma(2,3,11)$ are
contactomorphic.)
\smallskip\\
Furthermore, at least one of the contact structures $\xi_{\pm}$ on $+\Sigma(2,3,6n-1)$, $n\geq2$,
is not planar, and hence Corollary \ref{briecor} does not follow from the result in \cite{ach}.
Indeed, $M_c(2,3,6n-1)$ carries a Stein structure (see \cite[Section 1]{milnorfill}), which induces a
contact structure on $\Sigma(2,3,6n-1)$ and $\bp\big(M_c(2,3,6n-1)\big)=2(n-1)$ (see \cite[p.~74]{kirby}).
(An alternative strategy for $n=2$ can be found in \cite[Proof of Theorem 1.9]{geo}.)
With \cite[Theorem 4.1]{planar}, which tells us that $\bp\big(M_c(2,3,6n-1)\big)$ must vanish in the planar case,
the claim follows.
\smallskip\\
On $-\Sigma(2,3,11)$ there exists exactly one tight contact structure $\xi_0$, and the contact
structure $\xi_0$ is Stein fillable (see \cite[Theorem 4.9]{2311}).
It is not known whether there exists a Stein filling of $(-\Sigma(2,3,11),\xi_0)$ with non trivial Chern
class.
In \cite[Remark 3.3.3]{gauge} it is conjectured that any Stein filling of $(-\Sigma(2,3,11),\xi_0)$ is
diffeomorphic to the nucleus $N(2)$.
Notice that the only possible Stein structure on $N(2)$ has trivial first Chern class (use that $N(2)$ is
simply connected, \cite[Figure 12.81]{kirby} and Theorem \ref{egconstr}). 
Therefore we do not know whether our criterion applies in this situation or not.
Further, by \cite[Theorem 1.3(c) and Corollary 2.2]{transfol} $\xi_0$ is universally tight (use that 
$1/2 < 3/5$, $1/3 < 2/$5 and $2/11 < 1/5$) and not planar (use again \cite[Theorem 4.1]{planar}
and that $\bp\big(N(2)\big)=1$; see argumentation after Remark \ref{bpeq}).
\smallskip\\
For a discussion of $-\Sigma(2,3,17)$ we refer the reader to \cite{ozfill}.
Further we remark that our technique may apply to small Seifert fibred manifolds as studied
for example in \cite{gls,wu}.


\section{Non-contractible Reeb orbits on odd lens spaces \label{lens}}

For coprime natural numbers $p$ and $q$ satisfying $p>q\geq1$ the lens space $L_{p,q}$
is defined as the quotient $\S^3/G_{p,q}$ of the unit sphere in $\C^2$ by the discrete
subgroup $G_{p,q}$ of all diagonal matrices $\diag(\zeta,\zeta^q)\in\uni(2)$ with $\zeta^p=1$.
Notice that $\pi_2(L_{p,q})=0$ and $\pi_1(L_{p,q})=\Z_p=H_1(L_{p,q};\Z)$.
There exists a unique continued fraction expansion
                                   \begin{gather*}
                                     -\frac{p}{q}=
                                     n_1-\cfrac{1}{n_2-\cfrac{1}{\dotsb -\cfrac{1}{n_k} }}
                                     \;,
                                     \qquad
                                     \mtq{with integers}
                                     n_i<-1
                                     \;,\;
                                     i\in\{1,\ldots,k\}
                                     \;,
                                   \end{gather*}
of $-\frac{p}{q}$ which will be shortly denoted by $[n_1,\ldots,n_k]$.
The classification theorem due to Honda \cite{honda1} for tight contact structures on lens spaces
$L_{p,q}$ states that there exist exactly
                                  \begin{gather*}
                                    |(n_1+1)(n_2+1)\cdots (n_k+1)|
                                  \end{gather*}
tight contact structures up to isotopy.
All the tight contact structures on $L_{p,q}$ are obtained by Legendrian $(-1)$-surgery on
Legendrian links $L$ in $\S^3$ (see \cite{honda1}).
These links $L$ are linked chains $\big((K_1,n_1),\ldots,(K_k,n_k)\big)$ of framed unknots
(\ie they admit a Seifert surface diffeomorphic to $\D^2$) with $n_i=\tb(K_i)-1$ for all
$i\in\{1,\ldots,k\}$ as shown in \cite[Figure 16]{honda1}.
Therefore, all tight contact structures on $L_{p,q}$ are Stein fillable by Theorem \ref{egconstr}.
The rotation number $\rot(K_i)$, for each $i=1,\ldots,k$, can have any of the following values
                                    \begin{gather*}
                                      n_i+2\;,\; n_i+4\;,\;\ldots\;,\; n_i+2|n_i+1|
                                     \;,
                                    \end{gather*}
which are exactly the values allowed by the Bennequin-inequality
$\tb(K_i)+|\rot(K_i)|\leq -1$ and the condition $\tb(K_i)+\rot(K_i)\equiv 1\pmod 2$,
\cf \cite{legtranstifgtel}.
\smallskip\\
Recall, that the lens space $L_{p,q}=L_{p,q}^{\odd}$ is called odd if there exists at least one odd
integer $n_i$ in the associated continued fraction expansion $[n_1,\ldots,n_k]$.
As we already remarked in the introduction, the Weinstein conjecture holds true for all lens spaces.
Alternatively, for odd lens spaces the Weinstein conjecture follows from Corollary \ref{wcegconstr} with the
above mentioned isotopy classification theorem (notice that if the Weinstein conjecture holds true for $\xi$
then the Weinstein conjecture follows trivially for all contact structures contactomorphic to $\xi$) and the
verified Weinstein conjecture for overtwisted contact structures (see \cite{h2}).
The following result is not covered by Hofer's approach in \cite{h2}.

\begin{prop}
  \label{homol}
  For all odd lens spaces $L_{p,q}^{\odd}$ and all tight contact structures $\xi$ there exists a
  $\lambda$-Reeb link not homologous to zero in $H_1(L_{p,q}^{\odd};\Z)$ for all $\lambda\in\pc(\xi_+)$.
\end{prop}

\begin{proof}[\bf Proof of Theorem \ref{nclroools}]
  By Proposition \ref{homol} for all $\lambda\in\pc(\xi_+)$ there exists a $\lambda$-Reeb link not homologous
  to zero in $H_1(L_{p,q}^{\odd};\Z)$.
  In particular, there exists a component of the $\lambda$-Reeb link which is not contractible.
\end{proof}

\begin{proof}[\bf Proof of Proposition \ref{homol}]
  Let $W$ be the Stein filling manifold of $(L_{p,q}^{\odd},\xi_{r_1,\ldots,r_k})$ obtained from $\D^4$ via
  Legendrian $(-1)$-surgery along the framed link $\big((K_1,n_1),\ldots,(K_k,n_k)\big)$
  in $\S^3=\partial\D^4$ as described above, with
                                                            \begin{gather*}
                                                              r_j\define\rot(K_j)
                                                              \mtq{for all}
                                                              j\in\set{1,\ldots,k}
                                                              \;.
                                                            \end{gather*}
  Recall that $n_j\leq -2$ and notice that
                                                            \begin{gather}
                                                              \label{homol3}
                                                              n_j+r_j\equiv 0\pmod 2
                                                              \mtqq{and}
                                                              |r_j|\leq-n_j-2
                                                              \;.
                                                            \end{gather}
  We orient the knots $K_j$ in such a way that the linking matrix takes the form
                                                            \begin{gather*}
                                                              \big(\lk(K_i,K_j)\big)=
                                                              \left( 
                                                                \begin{array}{ccccc} 
                                                                  n_1 & 1          \\
                                                                  1      & n_2     & 1  \\
                                                                  &      & \ddots \\
                                                                  &      &   1         & n_{k-1} & 1     \\
                                                                  &      &               &  1           & n_k 
                                                                \end{array} 
                                                              \right)
                                                            \end{gather*} 
  with respect to the (free) basis $[C_1],\ldots,[C_k]$ of $H_2(W,\partial W;\Z)$, \cf Remark \ref{compute}.
  The equations \eqref{homol3} as well as the equations \eqref{orient1}, \eqref{orient2} and \eqref{lkrel} in
  Remark \ref{compute} are valid regardless of which orientations we choose.
  \smallskip\\
  The image of $[C_j]$ in $H_1(L_{p,q}^{\odd};\Z)$ under the connecting homomorphism is denoted
  by $c_j$.
  Because at least one of the $n_j$'s is odd the corresponding $r_j$ does not vanish.
  Therefore, by Corollary \ref{wcegconstr}, there exists a $\lambda$-Reeb link in $L_{p,q}^{\odd}$
  representing the integral $1$-homology class
                                                            \begin{gather*}
                                                              x\define
                                                              -\sum_{j=1}^kr_jc_j
                                                              \;,
                                                            \end{gather*}
  for all $\lambda\in\pc(\xi_{r_1,\ldots,r_k})$.
  Using the relations \eqref{lkrel} in Remark \ref{compute} and $c_0=0=c_{k+1}$ we get
                                                            \begin{gather*}
                                                              c_{j+1}=
                                                              -c_{j-1}-n_jc_j
                                                              \;,
                                                              \mtq{for}
                                                              j=1,\ldots,k
                                                              \;.
                                                            \end{gather*}
  Hence, there are integers $e_j$, $j=1,\ldots,k$, unique modulo $p$ such that $c_{k+1-j}=e_jc_k$.
  With $e_0=0=e_{k+1}$ we get
                                                            \begin{gather}
                                                              \label{homol1}
                                                              e_1=1
                                                              \mtqq{and}
                                                              e_{j+1}=
                                                              -e_{j-1}-n_{k+1-j}e_j
                                                              \;,
                                                              \mtq{for}
                                                              j=1,\ldots,k
                                                              \;,
                                                            \end{gather}
  as well as
                                                            \begin{gather*}
                                                              x=
                                                              \bigg(
                                                                \sum_{j=1}^kr_je_{k+1-j}
                                                              \bigg)
                                                              (-c_k)
                                                              \in H_1(L_{p,q}^{\odd};\Z)
                                                              \;.
                                                            \end{gather*}
  Further, there are uniquely determined coprime integers $p_j>q_j\geq1$ defined by
                                                            \begin{gather*}
                                                              -\frac{p_j}{q_j}=
                                                              [n_{k+1-j},\ldots,n_k]
                                                              \mtq{for}
                                                              j=1,\ldots,k
                                                              \;.
                                                            \end{gather*}
  Notice that with $q_0=0$ we get
                                                            \begin{gather}
                                                              \label{homol2}
                                                              q_1=1
                                                              \mtqq{and}
                                                              q_{j+1}=
                                                              -q_{j-1}-n_{k+1-j}q_j
                                                              \mtq{for}
                                                              j=1,\ldots,k
                                                              \;,
                                                            \end{gather}
  as well as $q_k=q$, $q_{k+1}=p$ and $q_{j+1}>q_j\geq 1$ for all $j=1,\ldots,k$.
  By \eqref{homol1} and \eqref{homol2} we find
                                                            \begin{gather*}
                                                              e_j\equiv q_j \pmod p
                                                              \;,
                                                              \mtq{for}
                                                              j=0,\ldots,k+1
                                                              \;.
                                                            \end{gather*}
  The claim is equivalent to $x\not=0$ in $H_1(L_{p,q}^{\odd};\Z)$ or
  $\sum_{j=1}^kr_je_{k+1-j}\not\equiv 0\pmod p$.
  Arguing by contradiction we suppose that $x=0$.
  Then, representing the residual classes $e_j$ by the integers $q_j$ for all $j=0,\ldots,k+1$, either
                                                            \begin{gather}
                                                              \label{alternative}
                                                              \sum_{j=1}^kr_jq_{k+1-j}=0
                                                              \mtqq{or}
                                                              p\leq
                                                              \bigg|
                                                                \sum_{j=1}^kr_jq_{k+1-j}
                                                              \bigg|
                                                              \;.
                                                            \end{gather}
  Supposing that the latter is true we get, using \eqref{homol3} and \eqref{homol2},
                                                            \begin{eqnarray*}
                                                              p
                                                              &\leq&
                                                              \sum_{j=1}^k|r_j|q_{k+1-j}\leq
                                                              \sum_{j=1}^k(-n_j-2)q_{k+1-j}
                                                              \\
                                                              &=&
                                                              -2\sum_{j=1}^kq_{k+1-j}+
                                                              \sum_{j=1}^kq_{k-j}+
                                                              \sum_{j=1}^kq_{k+2-j}
                                                              \\
                                                              &=&
                                                              -2\sum_{j=1}^kq_{k+1-j}+
                                                              \sum_{j=2}^{k+1}q_{k+1-j}+
                                                              \sum_{j=0}^{k-1}q_{k+1-j}
                                                              \\
                                                              &=&
                                                              -2q_k-2q_1+q_1+q_0+q_{k+1}+q_k
                                                              \\
                                                              &=&
                                                              -q-1+p
                                                              \;.
                                                            \end{eqnarray*}
  This leads to $q\leq-1$ which is a contradiction.
  Therefore, with the first equation in \eqref{alternative} we get
                                                            \begin{gather}
                                                              \label{homol4}
                                                              -r_1q_k=
                                                              \sum_{j=2}^kr_jq_{k+1-j}
                                                              \;,
                                                            \end{gather}
  and so $q_k$ divides the right hand side of \eqref{homol4}.
  Then either
                                                            \begin{gather}
                                                              \label{homol5}
                                                              \sum_{j=2}^kr_jq_{k+1-j}=0
                                                              \mtqq{or}
                                                              q_k\leq
                                                              \bigg|
                                                              \sum_{j=2}^kr_jq_{k+1-j}
                                                              \bigg|
                                                              \;.
                                                            \end{gather}
  Supposing that the latter is true we get with a similar reasoning as above
                                                            \begin{eqnarray*}
                                                              q_k
                                                              &\leq&
                                                              -2\sum_{j=2}^kq_{k+1-j}+
                                                              \sum_{j=2}^kq_{k-j}+
                                                              \sum_{j=2}^kq_{k+2-j}
                                                              \\
                                                              &=&
                                                              -2q_{k-1}-2q_1+q_1+q_0+q_k+q_{k-1}
                                                              \\
                                                              &=&
                                                              -q_{k-1}-1+q_k
                                                              \;.
                                                            \end{eqnarray*}
  This leads to $q_{k-1}\leq-1$ which is a contradiction.
  So the first case in \eqref{homol5} is left.
  Therefore,  \eqref{homol4} gives $r_1q_k=0$, \ie $r_1=0$ and hence, by \eqref{homol3}, $n_1\equiv 0\pmod 2$.
  If we repeat this argument we end up with $r_k=0$ and hence $n_k\equiv 0\pmod 2$.
  This shows that all $n_j$'s are even, contradicting our assumption on $L_{p,q}^{\odd}$ to be odd.
\end{proof}


\section{A directed cobordism \label{dircob}}

A hypersurface $M$ of a connected manifold $W$ is called {\it separating}\/ if $W\setminus M$
is not connected.
In the case of a separating hypersurface of contact type Theorem \ref{base} can be slightly
extended to the following

\begin{prop}
  \label{corchen}
  Let $(W,\omega)$ be a closed connected symplectic $4$-manifold satisfying $\bp(W)>1$
  and let $(M,\xi_+)$ be a separating hypersurface of contact type in $(W,\omega)$.
  Then the implications \eqref{base1} and \eqref{base2} of Theorem \ref{base} hold for any
  $\xi_+$-defining contact form $\lambda\in\pc(\xi_+)$.
\end{prop}

We shall prove Proposition \ref{corchen} by a gluing argument used before in \cite{McCW,elunique,es}
and by directed cobordisms.
Let $(W,\omega)$ be a compact symplectic $4$-manifold with boundary
$\partial W=\overline{M}_1\cup M_2$ and $M_1,M_2\neq\emptyset$.
If $(M_1,\lambda_1)$ and $(M_2,\lambda_2)$ are hypersurfaces of contact type, then
$(W,\omega)$ is called a {\it directed symplectic cobordism}\/ from $(M_1,\lambda_1)$
to $(M_2,\lambda_2)$ and we will write $(M_1,\lambda_1)\prec_{\omega} (M_2,\lambda_2)$
instead of $(W,\omega)$ in this case.
Notice that our terminology is borrowed from \cite{ehdscob} and differs from the one
used in \cite{egh}.
Again, if the $\xi_+$-defining contact form $\lambda$ is not needed, we will use $\xi_+$
in our notation.
Similar to the case of Stein fillings we call a directed symplectic cobordism
$(M_1,\lambda_1)\prec_{\omega} (M_2,\lambda_2)$ a {\it Stein cobordism}\/ if it carries a Stein
structure.

\begin{lem}
  \label{usticob}
  Let $(M,\xi_+)$ be a closed contact manifold and $\lambda_{j}\in\pc(\xi_+)$, $j=1,2$,
  $\xi_+$-defining contact forms.
  Then there exists a positive constant $c_{12}>0$ and a directed symplectic cobordism
  $(M,\lambda_1)\prec_{\omega_{12}}(M,c_{12}\lambda_2)$ diffeomorphic to $[-1,1]\times M$.
\end{lem}

\begin{proof}
  We use a construction from Ustilovsky's thesis \cite[Section 3.6]{ustia}.
  We consider the {\it symplectisation}\/ $\big(\R\times M,d(e^{\theta}\lambda_1)\big)$ of $(M,\lambda_1)$.
  There exists a function $f_{12}\in C^{\infty}(M,\R_{>0})$ on the compact manifold $M$ such that
  $\lambda_1=f_{12}\lambda_2$.
  Let $R>0$ be a constant which will be chosen later and $\beta\in C^{\infty}(\R,[0,1])$
  such that $\beta|_{(-\infty,-1]}=0$, $\beta|_{[1,\infty)}=1$ and $\beta'\geq 0$.
  Define a function $f\in C^{\infty}(\R\times M,\R_{>0})$ by
                                  \begin{gather*}
                                    f(\theta,p)\define
                                    e^{\theta+R}
                                    \Big(
                                    \big(
                                    1-\beta(\tfrac{\theta}{R})
                                    \big)
                                    f_{12}(p)+\beta(\tfrac{\theta}{R})
                                    \Big)\;,
                                    \qquad
                                    (\theta,p)\in\R\times M\;.
                                  \end{gather*}
  Note that
                                 \begin{gather*}
                                   \partial_{\theta}f(\theta,p)=
                                   e^{\theta+R}
                                   \Big(
                                   \tfrac{1}{R}\beta'(\tfrac{\theta}{R})
                                   \big(
                                   1-f_{12}(p)
                                   \big)+
                                   \big(
                                   1-\beta(\tfrac{\theta}{R})
                                   \big)
                                   f_{12}(p)+\beta(\tfrac{\theta}{R})
                                   \Big)
                                 \end{gather*}
  for all $(\theta,p)\in\R\times M$.
  There exists $R_{12}>0$ such that for all $R\geq R_{12}$ we have $\partial_{\theta}f>0$ on
  $\R\times M$.
  Consider the closed $2$-form $\omega_{12}\define d(f\lambda_2)$ on $\R\times M$ with $f$
  defined using $R=R_{12}$.
  Then
                                \begin{gather*}
                                  \omega_{12}\wedge\omega_{12}=
                                  2fdf\wedge\lambda_2\wedge d\lambda_2
                                  \qquad
                                  \text{and hence}
                                  \qquad
                                  \iota_{\partial_{\theta}}(\omega_{12}\wedge\omega_{12})=
                                  2f(\partial_{\theta}f)\lambda_2\wedge d\lambda_2\;.
                                \end{gather*}
  Therefore, $\omega_{12}$ is a symplectic form on $\R\times M$ which equals 
  $d(e^{\theta+R_{12}}\lambda_1)$ on $(-\infty,-R_{12}]\times M$ and
  $d(e^{\theta+R_{12}}\lambda_2)$ on $[R_{12},\infty)\times M$.
  The symplectic manifold $([-R_{12},R_{12}]\times M,\omega_{12})$ defines the claimed cobordism
  $(M,\lambda_1)\prec_{\omega_{12}} (M,c_{12}\lambda_2)$ with $c_{12}=e^{2R_{12}}$.
\end{proof}

\begin{proof}[\bf Proof of Proposition \ref{corchen}]
  Let $(M,\lambda_1)$ be a separating hypersurface of contact type in $(W,\omega)$ and
  $\lambda=\lambda_2\in\pc(\xi_+)$.
  Denote by $W_{\pm}$ the closures of the components of $W\setminus M$ and
  $\omega_{\pm}\define\omega\restrict_{W_{\pm}}$, where the sign is chosen such that 
  $(W_-,\omega_-)$ is the strong convex filling of $(M,\lambda_1)$ and $(W_+,\omega_+)$ is the
  strong concave filling of $(M,\lambda_1)$.
  There exist collar neighbourhoods $U_{\pm}$ of $M$ in $W_{\pm}$ such that in the notation of
  the proof of Lemma \ref{usticob} we have that $(U_-,c_{12}^{-1}\omega_-)$ is symplectomorphic to
  $\big((-\varepsilon,0]\times M,c_{12}^{-1}d(e^{\theta}\lambda_1)\big)$
  and $(U_+,c_{21}\omega_+)$ is symplectomorphic to
  $\big([0,\varepsilon)\times M,c_{21}d(e^{\theta}\lambda_1)\big)$
  for some $\varepsilon>0$ (\cf \cite[Section 2]{es}).
  Denote the corresponding symplectomorphisms by $\varphi_{\pm}$.
  Gluing $(W_-,c_{12}^{-1}\omega_-)$, 
  $(M,c_{12}^{-1}\lambda_1)\prec_{c_{12}^{-1}\omega_{12}} (M,\lambda_2)$,
  $(M,\lambda_2)\prec_{\omega_{21}} (M,c_{21}\lambda_1)$ and $(W_+,c_{21}\omega_+)$
  along the boundaries using the symplectomorphisms $\varphi_{\pm}$ yields a symplectic
  manifold $(W',\omega')$ such that $(M,\lambda_2)$ is a hypersurface of contact type in
  $(W',\omega')$.
  Proposition \ref{corchen} follows now from Theorem \ref{base} because $W$ and $W'$ are
  homotopy equivalent.
\end{proof}

The construction in the proof of Proposition \ref{corchen} can be used to glue directed
symplectic cobordisms or both kinds of strong fillings along orientation-reversed
contactomorphic contact manifolds appearing as boundary components.
For that one must allow rescaling by positive constants of the corresponding symplectic or
contact forms.
Therefore, only the involved co-oriented contact structures and positive conformal classes
of symplectic forms are respected.
We recall the gluing construction from \cite{elunique} (which also follows from the proof of
Proposition \ref{corchen}):

\begin{defi}[Gluing along boundaries of contact type]
  Let $(W_j,\omega_j)$ be a symplectic manifold with nonempty boundary and let $(M_j,\xi_+^j)$
  be a hypersurface of contact type in $(W_j,\omega_j)$, $j=1,2$,
  such that $M_j$ is a boundary component of $W_j$.
  Suppose that there exists an orientation-reversing contactomorphism
  $\varphi_{12}:(M_1,\xi_+^1)\lra(M_2,\xi_+^2)$.
  The manifold $W_1\cup_{\varphi_{12}}W_2$ obtained by gluing along $M_1$ via $\varphi_{12}$
  carries a symplectic form whose restriction to $W_1$ coincides with $\omega_1$.
  The resulting symplectic manifold is denoted by $(W_1\cup_{\xi_+^1}W_2,\omega_{\xi_+^1})$.
\end{defi}


\section{Realisation and proof of the main theorems \label{rpmt}}

An embedding $f:(W_1,\omega_1)\longrightarrow (W_2,\omega_2)$ of a symplectic manifold
$(W_1,\omega_1)$ into a symplectic manifold $(W_2,\omega_2)$ is called {\it iso-symplectic}\/
if $f^*\omega_2=\omega_1$.

\begin{prop}
  \label{realbp}
  Let $(M,\xi_+)$ be a closed contact $3$-manifold and $(W,\omega)$ be a strong convex
  filling of it.
  Then for every $n\in\N$ there exists a closed connected symplectic $4$-manifold $\big(W(n),\omega(n)\big)$ 
  with $\bp\big(W(n)\big)\geq n$ such that $(W,\omega)$ admits an iso-symplectic embedding into
  $\big(W(n),\omega(n)\big)$.
  In particular $(M,\xi_+)$ is contactomorphic to a separating hypersurface of contact
  type in $\big(W(n),\omega(n)\big)$.
\end{prop}

\begin{proof}
  Set $(M_0,\xi_+^0)=(M,\xi_+)$ and $(W_0,\omega_0)=(W,\omega)$.
  Let $(M_j,\xi_+^j)$, $j\in\{1,\ldots,n\}$, be a finite collection of closed connected contact manifolds,
  which we specify below.
  The iterated connected sum $M_{\#^n}=\#_{j=0}^nM_j$ carries a co-oriented contact structure
  denoted by $\xi_+^{\#_n}$.
  There exists a directed symplectic cobordism
                                \begin{gather*}
                                  (W_{\#^n},\omega_{\#^n})\define
                                  (M_{\sqcup^n},\xi_+^{\sqcup_n})
                                  \prec
                                  (M_{\#^n},\xi_+^{\#_n})\;,
                                \end{gather*}
  where $(M_{\sqcup^n},\xi_+^{\sqcup_n})=\bigsqcup_{j=0}^n(M_j,\xi_+^j)$ (see \cite{wst1}).
  (If $M_0$ is not connected, perform additional connected sum surgeries along the ordered
  components.
  Hence $(M_0,\xi_+^0)$ is directed cobordant to a connected contact manifold.
  For simplicity we assume that $M_0$ is already connected.)
  Let $(W_{\cav},\omega_{\cav})$ be a strong concave filling of $(M_{\#^n},\xi_+^{\#_n})$
  ensured by Theorem \ref{concave} and set
                              \begin{gather*}
                                (\widehat{W},\widehat{\omega})\define
                                W_0
                                \cup_{\xi_+^0}
                                W_{\#^n}
                                \cup_{\xi_+^{\#_n}}
                                W_{\cav}\;.
                              \end{gather*}
  Then $\widehat{W}$ is connected by construction.
  Suppose that for each $j\in\{1,\ldots,n\}$ there exists a connected strong convex filling
  $(W_j,\omega_j)$ of $(M_j,\xi_+^j)$.
  Set
                          \begin{gather*}
                            \big(
                            W(n),\omega(n)
                            \big)\define
                            \textstyle{\bigsqcup_{j=1}^n}W_j
                            \cup_{\xi_+^{\sqcup_n}}
                            \widehat{W}\;.
                          \end{gather*}

  \begin{rem}
    \label{bpeq}
    A closed orientable connected $3$-manifold $Y$ is called an {\it (integral) homology sphere}\/ if the
    first integral homology vanishes, \ie $H_1(Y;\Z)=H_1(\S^3;\Z)$.
    If $M_j$ is a homology sphere for each $j\in\{1,\ldots,n\}$ then
                         \begin{gather*}
                           \bp\big(W(n)\big)=
                           \bp(\widehat{W})+\sum_{j=1}^n\bp(W_j)\;,
                         \end{gather*} 
    because the intersection form splits as
    $Q_{W(n)}=Q_{\widehat{W}}\oplus Q_{W_1}\oplus\ldots\oplus Q_{W_n}$
    (see \cite[Exercise 1.3.5.(b)$^*$]{kirby}).
  \end{rem}

  For example we take the homology sphere $M_j\define\partial N(2)$; the boundary of the
  Gompf nucleus $N(2)$ (see \cite{nuclei}).
  $N(2)$ carries the structure of a Stein manifold with boundary inducing a contact structure
  on $\partial N(2)$ and satisfies $\bp\big(N(2)\big)=1$
  (see \cite[Theorem 1.7 and Remark 3.3.1]{gauge} or \cite[p.~515]{tightsw}).
  Proposition \ref{realbp} follows from Remark \ref{bpeq}.
\end{proof}

\begin{proof}[\bf Proof of Theorem \ref{Chen1} and Theorem \ref{Chen2}]
  The claims follow from Proposition \ref{realbp} (with $n=2$) and Proposition \ref{corchen}.
\end{proof}

We remark that it is known (in principle) how to construct symplectic $4$-manifolds with $\bp>1$
via surgery.
For example the proof of the second part of Theorem \ref{concave} includes a method via fibre sum
with elliptic surfaces along symplectic tori.
In \cite{wittens,genus} similar statements can be found.
\smallskip\\
We sketch a further approach.
Proposition \ref{realbp} can be obtained by the following argument whose germ was already used in
\cite[Theorem 3.2]{tightsw} (see \cite[Lemma 3.1]{etnfll} and \cite[Lemma 3.1]{geo}).
If $M$ is not a homology sphere then glue a directed cobordism $(M,\xi_+)\prec (M_1,\xi^1_+)$
to $(W,\omega)$ along $(M,\xi_+)$ where $M_1$ is a homology sphere depending on $M$
(we used the notation from Proposition \ref{realbp} and assumed for simplicity that $M$ is connected).
Using Legendrian $(-1)$-surgery one can find inductively directed cobordisms $(W_n,\omega_n)$
from $ (M_n,\xi^n_+)$ to $(M_{n+1},\xi^{n+1}_+)$ with $\bp(W_n)\geq 1$, $n\in\N$, where
$M_{n+1}$ is again a homology sphere.
Using a concave filling of $(M_{n+1},\xi^{n+1}_+)$, gluing along the contact type boundaries yields a 
closed connected symplectic $4$-manifold $\big(W(n),\omega(n)\big)$ with $\bp\big(W(n)\big)\geq n$ 
(see Remark \ref{bpeq}).

\def\cprime{$'$} \def\cprime{$'$} \def\cprime{$'$} \def\cprime{$'$}


\begin{thebibliography}{10}

\bibitem{ach}
{\sc C.~Abbas, K.~Cieliebak, and H.~Hofer}, {\em {The Weinstein Conjecture for
  Planar Contact Structures in Dimension Three}}.
\newblock arXiv:math.SG/0409355.

\bibitem{aostein}
{\sc S.~Akbulut and B.~\"Ozba\u{g}c\i}, {\em On the topology of compact {S}tein
  surfaces}, Int. Math. Res. Not.,  (2002), pp.~769--782.

\bibitem{botttu}
{\sc R.~Bott and L.~W. Tu}, {\em Differential forms in algebraic topology},
  vol.~82 of Graduate Texts in Mathematics, Springer-Verlag, New York, 1982.

\bibitem{milnorfill}
{\sc C.~Caubel, A.~Nemethi, and P.~Popescu-Pampu}, {\em {Milnor open books and
  Milnor fillable contact 3-manifolds}}.
\newblock arXiv:math.SG/0409160.

\bibitem{weichen}
{\sc W.~Chen}, {\em Pseudo-holomorphic curves and the {W}einstein conjecture},
  Comm. Anal. Geom., 8 (2000), pp.~115--131.

\bibitem{dgs}
{\sc F.~Ding, H.~Geiges, and A.~I. Stipsicz}, {\em Surgery diagrams for contact
  3-manifolds}, Turkish J. Math., 28 (2004), pp.~41--74.

\bibitem{elfll}
{\sc Y.~Eliashberg}, {\em Filling by holomorphic discs and its applications},
  in Geometry of low-dimensional manifolds, 2 (Durham, 1989), vol.~151 of
  London Math. Soc. Lecture Note Ser., Cambridge Univ. Press, Cambridge, 1990,
  pp.~45--67.

\bibitem{stein}
\leavevmode\vrule height 2pt depth -1.6pt width 23pt, {\em Topological
  characterization of {S}tein manifolds of dimension {$>2$}}, Internat. J.
  Math., 1 (1990), pp.~29--46.

\bibitem{el1}
\leavevmode\vrule height 2pt depth -1.6pt width 23pt, {\em Contact
  {$3$}-manifolds twenty years since {J}. {M}artinet's work}, Ann. Inst.
  Fourier (Grenoble), 42 (1992), pp.~165--192.

\bibitem{legtranstifgtel}
\leavevmode\vrule height 2pt depth -1.6pt width 23pt, {\em Legendrian and
  transversal knots in tight contact {$3$}-manifolds}, in Topological methods
  in modern mathematics (Stony Brook, NY, 1991), Publish or Perish, Houston,
  TX, 1993, pp.~171--193.

\bibitem{elunique}
\leavevmode\vrule height 2pt depth -1.6pt width 23pt, {\em Unique
  holomorphically fillable contact structure on the {$3$}-torus}, Internat.
  Math. Res. Notices,  (1996), pp.~77--82.

\bibitem{few}
\leavevmode\vrule height 2pt depth -1.6pt width 23pt, {\em A few remarks about
  symplectic filling}, Geom. Topol., 8 (2004), pp.~277--293 (electronic).

\bibitem{egh}
{\sc Y.~Eliashberg, A.~Givental, and H.~Hofer}, {\em Introduction to symplectic
  field theory}, Geom. Funct. Anal.,  (2000), pp.~560--673.
\newblock GAFA 2000 (Tel Aviv, 1999).

\bibitem{planarlec}
{\sc J.~B. Etnyre}, {\em {Lectures on open book decompositions and contact
  structures}}.
\newblock arXiv:math.SG/0409402.

\bibitem{planar}
\leavevmode\vrule height 2pt depth -1.6pt width 23pt, {\em {Planar open book
  decompositions and contact structures}}.
\newblock arXiv:math.SG/0404267.

\bibitem{es}
\leavevmode\vrule height 2pt depth -1.6pt width 23pt, {\em Symplectic convexity
  in low-dimensional topology}, Topology Appl., 88 (1998), pp.~3--25.
\newblock Symplectic, contact and low-dimensional topology (Athens, GA, 1996).

\bibitem{etnfll}
\leavevmode\vrule height 2pt depth -1.6pt width 23pt, {\em On symplectic
  fillings}, Algebr. Geom. Topol., 4 (2004), pp.~73--80 (electronic).

\bibitem{ehdscob}
{\sc J.~B. Etnyre and K.~Honda}, {\em On symplectic cobordisms}, Math. Ann.,
  323 (2002), pp.~31--39.

\bibitem{schlfrau}
{\sc U.~Frauenfelder and F.~Schlenk}, {\em {Applications of Hofer's geometry to
  Hamiltonian dynamics}}.
\newblock arXiv:math.SG/0305146.

\bibitem{explicit}
{\sc D.~T. Gay}, {\em Explicit concave fillings of contact three-manifolds},
  Math. Proc. Cambridge Philos. Soc., 133 (2002), pp.~431--441.

\bibitem{exageig}
{\sc H.~Geiges}, {\em Examples of symplectic {$4$}-manifolds with disconnected
  boundary of contact type}, Bull. London Math. Soc., 27 (1995), pp.~278--280.

\bibitem{ozfill}
{\sc P.~Ghiggini}, {\em {Ozsvath-Szabo invariants and fillability of contact
  structures}}.
\newblock arXiv:math.GT/0403367.

\bibitem{gls}
{\sc P.~Ghiggini, P.~Lisca, and A.~I. Stipsicz}, {\em {Classification of tight
  contact structures on small Seifert $3$-manifolds with $e_0\geq 0$}}.
\newblock arXiv:math.SG/0406080.

\bibitem{2311}
{\sc P.~Ghiggini and S.~Sch{\"o}nenberger}, {\em On the classification of tight
  contact structures}, in Topology and geometry of manifolds (Athens, GA,
  2001), vol.~71 of Proc. Sympos. Pure Math., Amer. Math. Soc., Providence, RI,
  2003, pp.~121--151.

\bibitem{ginz}
{\sc V.~L. Ginzburg}, {\em {The Weinstein conjecture and the theorems of nearby
  and almost existence}}.
\newblock arXiv:math.DG/0310330.

\bibitem{nuclei}
{\sc R.~E. Gompf}, {\em Nuclei of elliptic surfaces}, Topology, 30 (1991),
  pp.~479--511.

\bibitem{gompf}
\leavevmode\vrule height 2pt depth -1.6pt width 23pt, {\em Handlebody
  construction of {S}tein surfaces}, Ann. of Math. (2), 148 (1998),
  pp.~619--693.

\bibitem{kirby}
{\sc R.~E. Gompf and A.~I. Stipsicz}, {\em {$4$}-manifolds and {K}irby
  calculus}, vol.~20 of Graduate Studies in Mathematics, American Mathematical
  Society, Providence, RI, 1999.

\bibitem{3m}
{\sc A.~Hatcher}, {\em Notes on basic $3$-manifold topology}.
\newblock www.math.cornell.edu/${\sim}$hatcher/, 3M.pdf.

\bibitem{h2}
{\sc H.~Hofer}, {\em Pseudoholomorphic curves in symplectizations with
  applications to the {W}einstein conjecture in dimension three}, Invent.
  Math., 114 (1993), pp.~515--563.

\bibitem{s3tight}
{\sc H.~Hofer, K.~Wysocki, and E.~Zehnder}, {\em Finite energy foliations of
  tight three-spheres and {H}amiltonian dynamics}, Ann. of Math. (2), 157
  (2003), pp.~125--255.

\bibitem{hz}
{\sc H.~Hofer and E.~Zehnder}, {\em Symplectic invariants and {H}amiltonian
  dynamics}, Birkh\"auser Advanced Texts: Basler Lehrb\"ucher. [Birkh\"auser
  Advanced Texts: Basel Textbooks], Birkh\"auser Verlag, Basel, 1994.

\bibitem{honda1}
{\sc K.~Honda}, {\em On the classification of tight contact structures. {I}},
  Geom. Topol., 4 (2000), pp.~309--368 (electronic).

\bibitem{wittens}
{\sc P.~B. Kronheimer and T.~S. Mrowka}, {\em Witten's conjecture and property
  {P}}, Geom. Topol., 8 (2004), pp.~295--310 (electronic).

\bibitem{tightsw}
{\sc P.~Lisca and G.~Mati{\'c}}, {\em Tight contact structures and
  {S}eiberg-{W}itten invariants}, Invent. Math., 129 (1997), pp.~509--525.

\bibitem{transfol}
\leavevmode\vrule height 2pt depth -1.6pt width 23pt, {\em Transverse contact
  structures on {S}eifert 3-manifolds}, Algebr. Geom. Topol., 4 (2004),
  pp.~1125--1144 (electronic).

\bibitem{McCW}
{\sc J.~D. McCarthy and J.~G. Wolfson}, {\em Symplectic gluing along
  hypersurfaces and resolution of isolated orbifold singularities}, Invent.
  Math., 119 (1995), pp.~129--154.

\bibitem{mcduconty}
{\sc D.~McDuff}, {\em Symplectic manifolds with contact type boundaries},
  Invent. Math., 103 (1991), pp.~651--671.

\bibitem{latt2}
{\sc L.~I. Nicolaescu}, {\em Lattice points inside rational simplices and the
  {C}asson invariant of {B}rieskorn spheres}, Geom. Dedicata, 88 (2001),
  pp.~37--53.

\bibitem{genus}
{\sc P.~Ozsv{\'a}th and Z.~Szab{\'o}}, {\em Holomorphic disks and genus
  bounds}, Geom. Topol., 8 (2004), pp.~311--334 (electronic).

\bibitem{rab}
{\sc P.~H. Rabinowitz}, {\em Periodic solutions of a {H}amiltonian system on a
  prescribed energy surface}, J. Differential Equations, 33 (1979),
  pp.~336--352.

\bibitem{iohspheres}
{\sc N.~Saveliev}, {\em Invariants for homology {$3$}-spheres}, vol.~140 of
  Encyclopaedia of Mathematical Sciences, Springer-Verlag, Berlin, 2002.
\newblock Low-Dimensional Topology, 1.

\bibitem{gauge}
{\sc A.~I. Stipsicz}, {\em Gauge theory and {S}tein fillings of certain
  3-manifolds}, Turkish J. Math., 26 (2002), pp.~115--130.

\bibitem{geo}
\leavevmode\vrule height 2pt depth -1.6pt width 23pt, {\em On the geography of
  {S}tein fillings of certain 3-manifolds}, Michigan Math. J., 51 (2003),
  pp.~327--337.

\bibitem{existence}
{\sc C.~H. Taubes}, {\em The {S}eiberg-{W}itten and {G}romov invariants}, Math.
  Res. Lett., 2 (1995), pp.~221--238.

\bibitem{ustia}
{\sc I.~Ustilovsky}, {\em Contact homology and contact structures on {$S\sp
  {4m+1}$}}, PhD thesis, Stanford University, August 1999.

\bibitem{wstrab}
{\sc A.~Weinstein}, {\em On the hypotheses of {R}abinowitz' periodic orbit
  theorems}, J. Differential Equations, 33 (1979), pp.~353--358.

\bibitem{wst1}
\leavevmode\vrule height 2pt depth -1.6pt width 23pt, {\em Contact surgery and
  symplectic handlebodies}, Hokkaido Math. J., 20 (1991), pp.~241--251.

\bibitem{wu}
{\sc H.~Wu}, {\em {Tight Contact Small Seifert Spaces with $e_0\neq0,-1,-2$}}.
\newblock arXiv:math.GT/0402167.

\end{thebibliography}
\end{document}